\documentclass{article}
 
\usepackage{amsmath}
\usepackage{amssymb}
\usepackage{amscd}
\usepackage[all]{xy}

\usepackage{abstract}

          \newtheorem{df}{Definition}[section]
          \newtheorem{pr}[df]{Proposition}
          \newtheorem{theorem}[df]{Theorem}
          \newtheorem{lem}[df]{Lemma}
          
          \newtheorem{rem}[df]{\it Remark}

          \newtheorem{prob}[df]{Problem}
          \newtheorem{exam}[df]{Example}

          \newcommand{\qed}{~$\Box$\newline\hskip 0.18cm}
          \newcommand{\proof}{
                             \noindent{\bf Proof.~~}}

            \makeatletter
            \@addtoreset{equation}{section}
            \makeatother
         
         \newcommand{\mapright}[1]{%
           \smash{\mathop{%
           \hbox to 1cm{\rightarrowfill}}\limits^{#1}}}
        
         \newcommand{\mapleft}[1]{%
           \smash{\mathop{%
           \hbox to 1cm{\leftarrowfill}}\limits_{#1}}}
           
         \newcommand{\shortmapright}[1]{%
           \smash{\mathop{%
           \hbox to .7cm{\rightarrowfill}}\limits^{#1}}}
        
         \newcommand{\shortmapleft}[1]{%
           \smash{\mathop{%
           \hbox to .7cm{\leftarrowfill}}\limits_{#1}}}
           
\begin{document}
 
       
       \title{Liftings of pseudo-reflection groups of 
       toric quotients of Krull schemes\thanks{This version is dated   July 23, 2017.}}        
         \author{{Haruhisa} {\sc Nakajima}\thanks{Partially supported by  
         Grant No. 26400019:  the Japan Society for the Promotion of Sciences.}\\
        \small { Department of Mathematics}, 
        \small {\sc  J. F. Oberlin University}\\
        \small {\sc Machida}, {\sc Tokyo}  194-0294,
         {\sc JAPAN}} 
         \date{ }
                 \maketitle

       
         \begin{abstract} Let $G$ be an affine algebraic group with 
 a reductive  identity component $G^{0}$ acting regularly on an affine Krull scheme
 $X = {\rm Spec} (R)$ over an algebraically closed field.   Let $T$ be an algebraic subtorus  of  $G$
 and suppose that ${\mathcal Q}(R)^{T}= {\mathcal Q}(R^{T})$
 of quotient fields. We will show: 
If $G$ is the centralizer of $T$ in $G$,  then the pseudo-reflections of the action of $G$
on $R^{T}$ can be lifted to those on $R$. This result is applied  to  partially generalize Chevalley-Serre and
Steinberg theorems on pseudo-reflection groups. 
 
 \medskip
            
         \noindent
         {\it MSC:} primary 13A50, 14R20, 20G05; secondary 14L30, 14M25

     \noindent  
    {\it Keywords:} algebraic torus; pseudo-reflection;
               Krull domain; ramification index; invariant theory
          \end{abstract}

{\small\tableofcontents}
    \normalsize

         \small \section{Introduction}\normalsize
         {\bf 1.A.} For a commutative ring $R$, let ${\mathcal Q} (R)$ denote the total quotient ring
         of $R$ and ${\rm Spec} (R)$ the affine scheme defined by $R$.  Consider an
         action of  a  group $G$   on $R$ as automorphisms. For a prime ideal ${\mathfrak P}$
         of $R$, let 
         $$\begin{array}{l}{\mathcal D}_{G}({\mathfrak P}) =\{\sigma \in G \mid \sigma({\mathfrak P})
         ={\mathfrak P}\}\\
         {\mathcal I}_{G}({\mathfrak P})=  \{\sigma \in G \mid \sigma(x) -x \in {\mathfrak P}
     ~(x \in R) \}\end{array},$$
     which are respectively  referred to as the {\it decomposition group} and the {\it inertia group} of ${\mathfrak P}$
     under this action (cf. \cite{LR}).  For  homomorphisms $A \to B \to C$ of commutative rings, 
        let 
        ${\rm Ht}^{1} (C, A; B)$ be  defined to be   
        \begin{equation*}\left\{{\mathfrak p}\in {\rm Spec} (C) \mid {\rm ht} ({\mathfrak p})
=  {\rm ht} ({\mathfrak p}\cap B) =  {\rm ht} ({\mathfrak p}\cap A) =1\right\}.\end{equation*}
Here ${\rm ht} ({\mathfrak I})$ of an ideal ${\mathfrak I}$ stands for its height. 
In   case of  $B = C$, we   set 
  ${\rm Ht}^{1} (B, A):={\rm Ht}^{1} (C, A; B)$. 

         \vskip 0.15cm
 \noindent {\bf 1.B.}  In this paper algebraic groups are affine and defined  over  a fixed algebraically closed field $K$ of  an arbitrary characteristic $p$.     Affine $K$-schemes $X$ are 
        affine schemes of commutative $K$-algebras  $R$ 
         which are not necessarily   finite generated as algebras over $K$. 
           We say an action 
        $(X, G)$ or $(R, G)$ of  an affine algebraic group $G$ on  $X$
         is {\it regular}, when $G$ acts {\it rationally}
        on the $K$-algebra $R$ as $K$-algebra automorphisms (e.g., \cite{Po}).
         Furthermore
        $(X, G)$ is said to be {\it effective} if ${\rm Ker} (G \to {\rm Aut} (R))$ is finite.  
        If  a subset $S$ of $R$ is invariant under the action of $G$, 
        we denote by $G\vert_{S}$ the group consisting of the restriction 
        $\sigma\vert_{S}$ of all $\sigma\in G$ to $S$, which is called {\it  the
        group $G$ on $S$}. 
            
                \vskip 0.15cm
      \noindent {\bf 1.C.}   
       An affine $K$-scheme $X = {\rm Spec}  (R)$ is said to be Krull, if $R$ is a Krull 
       $K$-domain. For a prime ideal ${\mathfrak P}\in{\rm Ht}^{1}(R, R^{G})$, 
    let ${\rm e} ({\mathfrak P}, {\mathfrak P}\cap R^{G})$ be the ramification index
    ${\rm v}_{R, {\mathfrak P}}({\mathfrak P}\cap R^{G} )$,  where 
    ${\rm v}_{R, {\mathfrak P}}$ denotes the discrete valuation defined by 
    $R_{{\mathfrak P}}$. If   ${\mathfrak p}$ is a prime ideal of a Krull domain
    $R\cap L$ (e.g., \cite{BourbakiC}) for a subfield $L$ of ${\mathcal Q}(R)$ containing $K$  of ${\rm ht}({\mathfrak p}) =1$, 
    let $${\rm Over}_{{\mathfrak p}}(R) := \{ {\mathfrak P}\in {\rm Ht}^{1}(R, R\cap L) \mid
    {\mathfrak P}\cap (R\cap L) = {\mathfrak p}\},$$
    which is non-empty (e.g., \cite{Magid}). 
   The elements of ${\mathcal I}_{G}({\mathfrak P})$
      are referred to as the pseudo-reflections at ${\mathfrak P}$ under the action
     of $(X, G)$ (cf. \cite{Nak5}). This is a generalization of classical pseudo-reflections in Chap. IV of \cite{Bourbaki}. 
     If $\Gamma$ is a subset of a group, let $\langle \Gamma \rangle$ denote
     the subgroup generated by $\Gamma$. 
     For a closed normal subgroup $H$ of $G$,  we set  the subgroups 
    $$ 
     {\mathfrak R}(R^{H}, G) :=\left\langle  \bigcup _{{\mathfrak p} \in {\rm Ht}^{1} (R^{H}, R^{G})}{\mathcal I}_{G} ({\mathfrak p}) 
\right\rangle \subseteq G $$
 $$ {\mathfrak R}(R, G; H) :=\left\langle  \bigcup _{{\mathfrak P} \in {\rm Ht}^{1} (R, R^{G}; R^{H})}{\mathcal I}_{G} ({\mathfrak P}) 
\right\rangle \subseteq G
.$$
In case of  $H = \{1\}$, 
${\mathfrak R}(R, G)= {\mathfrak R}(R^{H}, G)$, which is called 
{\it the pseudo-reflection group of the action} $(X, G)$ or $(R, G)$. 
This group is {\it finite on} $X$ for $G$ with a reductive $G^{0}$ (cf. \cite{Nak5}),
which characterizes the reductively of $G^{0}$.
The purpose of this paper is to study {\it the  problem on   the lifting of pseudo-reflections} as follows:
\begin{prob}   Can any element  $\sigma$ of 
the   pseudo-reflection group ${\mathfrak R}(R^{H},  G)$ 
of the action $(R^{H}, G)$ be lifted to an element $\tilde{\sigma}$
of the one ${\mathfrak R}(R,  G)$ of
 $(R, G)$, i.e., $\tilde{\sigma}\vert _{R^{H}}=\sigma$?  
\end{prob} 



          \vskip 0.15cm           
 \noindent {{\bf 1.D.}}  Let  ${\mathfrak X}(G)$ be the group of rational characters of $G$ expressed as
an additive group with zero. For a rational $G$-module  $M$, put
\begin{equation*}
{\mathfrak X}(G)^{M} : = \{ \chi \in  {\mathfrak X}(G) \mid  M_{\chi}\not=\{0\}\}
\end{equation*}
 Here $M_{\chi}$ stands for 
  $\{ x\in M \mid \sigma (x) =\chi (\sigma) x\}$ of {\it relative invariants} 
  of $G$ in $M$ relative to $\chi$. 
 For a morphism $\gamma : H \to G$ of groups, let $Z_G(H)$ denote the
centralizer $\{\sigma \in G \mid \sigma\gamma(\tau) = \gamma (\tau) \sigma
 ~(\tau \in H)\}$ of a subset 
$\gamma(H)$ in  $G$. The main result of this paper is as follows: 
\begin{theorem}\label{maintheorem}
Suppose that $G^0$ is reductive and let $T$ be 
 a connected closed subgroup of $G^{0}$ which is   an algebraic torus. 
Let $(X, G)$ be an effective  regular action of $G$ on an affine Krull 
$K$-scheme $X$ with $X = {\rm Spec} (R)$. Suppose that $G = Z_{G}(T)$
and ${\mathcal Q}(R)^{T}= {\mathcal Q}(R^{T})$.   Then we have  
$${\mathfrak R}(R, G; T)\vert _{R^{T}} = {\mathfrak R}(R^{T}, G)\vert _{R^{T}}.$$ 
\end{theorem}
 
Thus Problem 1.1 is solved affirmatively in the case where $H = T$
under the certain condition as above. However this problem is not true 
in the case where $H$ is neither a torus nor $G = Z_{G}(T)$ that  are  discussed in Sect. 4.  In the proof of the main theorem, Proposition 3.3 plays an essential role and is shown by Theorem 3.3 of \cite{Nak3}   inspired by R.P. Stanley's Theorem 2.3 in \cite{Stanley}.  Furthermore in Sect. 5 we apply  the main theorem to partially generalize the  classical Chevalley-Serre  and
Steinberg theorems on pseudo-reflection groups. 


\section{Preliminaries}
\begin{lem}\label{normalcase}  Let $G$ be a group acting on an integrally closed domain $A$ as ring automorphisms and $N$ a normal subgroup of $G$ of a finite index. 
 \begin{itemize}
\item[(i)]  Let $H$ be a finite normal subgroup of $G$  such that $H \subseteq N$.  Then ${\mathcal Q}(A^G) = {\mathcal Q}(A)^G$ if and only if
 ${\mathcal Q}(A^N) = {\mathcal Q}(A^H)^N$. 
 \item[(ii)] For any prime ideal ${\mathfrak P}$ of $A$  the morphism
$A^{G}_{{\mathfrak P}\cap A^{G}} \to 
(A^{ {\mathcal I}_{G}({\mathfrak P}\cap A^{N})})_{{\mathfrak P}\cap A^{ {\mathcal I}_{G}({\mathfrak P}\cap
A^{N})}}$
is unramified in the sense of {\rm \cite{LR}}.  
\end{itemize}
 \end{lem}
\proof  {\it (i)}:  Suppose that ${\mathcal Q}(A^{G}) = {\mathcal Q}(A)^{G}$.  Let
$B$ be the integral closure of $A^{G}$ in ${\mathcal Q}(A)^{N}$. Then  $B \subseteq A\cap {\mathcal Q} (A)^{N}$.  Moreover we have 
${\mathcal Q}(B) = {\mathcal Q}(A)^{N}$, as ${\mathcal Q}(A)^{N}$
is Galois over ${\mathcal Q}(A)^{G}$.  
Consequently  ${\mathcal Q}(A^{N})= {\mathcal Q}(A)^{N}$.  The remainder of the proof  is omitted. 

{\it (ii)}: 
 Exchanging $N$ with ${\rm Ker} (G \to {\rm Aut} (A^{N}))$, we may suppose that 
 the  action $(B, G/N)$ on 
 $B = A^{N}$  is faithful and  Galois.  Put ${\mathfrak p}
={\mathfrak P}\cap B$.  Clearly ${\mathcal I}_{G}({\mathfrak p})\supseteq N$ and ${\mathcal I}_{G}({\mathfrak p})/N
= {\mathcal I}_{G/N} ({\mathfrak p})$. Then by Chap. VI of \cite{LR} 
$$(B^{G/N})_{{\mathfrak p}\cap B^{G/N}} \to 
(B^{({\mathcal I}_{G/N}({\mathfrak p})})_{{\mathfrak p}\cap B^{{\mathcal I}_{G/N}({\mathfrak p})} }$$
is unramified, which shows the assertion. 
\qed
\begin{pr}\label{finitecase}
Let $A$ be an integrally closed domain and $G$ a finite subgroup of ${\rm Aut} (A)$. Let  
${\mathfrak P}$ be a prime ideal of $A$. Suppose that $A^{G}_{{\mathfrak P}\cap A^{G}}$
is noetherian.  Then 
\begin{itemize}
\item[(i)]  For a subgroup $L$ of $G$, the canonical morphism
$A^{G}_{{\mathfrak P}\cap A^{G}} \to A^{L}_{{\mathfrak P}\cap A^{L}}$
is \'etale if and only if $L \supseteq {\mathcal I}_{G}({\mathfrak P})$. 
\item[(ii)]  For a normal subgroup $N$ of $G$, we have 
${\mathcal I}_{G}({\mathfrak P}) \cdot N = {\mathcal I}_{G}({\mathfrak P}\cap A^{N}).$
\end{itemize}
\end{pr}
\proof Replacing $A$ with $A^{G}_{{\mathfrak P}\cap A^{G}}\otimes_{A^{G}} A$, 
we may assume that  $A$ is noetherian, because $A$
is the integral closure of $A^{G}$ in a finite separable (Galois) extension ${\mathcal Q}(A)$ of ${\mathcal Q}(A)^{G}$. 
Then the assertion  {\it (i)} follows from
Proposition 2.2 and Corollary 2.4 of  Expos\'e V of  \cite{Gr}. 

 By {\it (i)}
we see that
\begin{equation}\label{eqnfinite}A^{G}_{{\mathfrak P}\cap A^{G}} \to (A^{{\mathcal I}_{G}({\mathfrak P})\cdot N})
_{{\mathfrak P}\cap A^{{\mathcal I}_{G}({\mathfrak P})\cdot N} }
\end{equation}
is \'etale. Consider the Galois group action $(B, G/N)$ with $B = A^{N}$ and put ${\mathfrak p}
={\mathfrak P}\cap B$. 
Then
  (\ref{eqnfinite})  is expressed  as   the \'etale morphism
$$(B^{G/N})_{{\mathfrak p}\cap B^{G/N}} \to 
(B^{({\mathcal I}_{G}({\mathfrak P})\cdot N)/N})_{{\mathfrak p}\cap B^{({\mathcal I}_{G}({\mathfrak P})\cdot N)/N} }.$$
From this and  {\it (i)} we infer that ${\mathcal I}_{G/N}({\mathfrak p})
\subseteq ({\mathcal I}_{G}({\mathfrak P})\cdot N)/N$. 
Since ${\mathcal I}_{G}({\mathfrak P})\cdot N \subseteq {\mathcal I}_{G}({\mathfrak P}\cap B)$
and ${\mathcal I}_{G/N}({\mathfrak p}) = 
{\mathcal I}_{G}({\mathfrak P}\cap B)/N$, we must have
${\mathcal I}_{G}({\mathfrak P})\cdot N = {\mathcal I}_{G}({\mathfrak P}\cap B)$.  \qed 
\begin{lem}\label{liftpart} Let $G$ be  an algebraic group with 
 an algebraic torus $G^{0}$ and suppose that $G = Z_{G}(G^{0})$. Then for any finite
subgroup $H$  of $G$, there exists a finite subgroup $M$ of $G$
such that $M \supseteq H$ and $G = M \cdot G^{0}$. 
\end{lem} 
\proof  For a semi-simple element $\tau \in G$, there is a finite subgroup $D_{\tau}$
of $G$ such that $D_{\tau}\cdot G^{0} = \langle G^{0}\cup \{\tau\}\rangle,$ 
because $\langle G^{0}\cup \{\tau\}\rangle$ is diagonalizable (e.g., \cite{Po}). For  $\sigma\in G$, let $(\sigma)_{s}$ 
be the semi-simple 
part of $\sigma$.  Let $\{ \sigma_{i} \in G \mid 
1\leqq i \leqq n \}$  be a complete set of representatives of $G/G^{0}$.   We denote by 
$M$ the subgroup of $G$ generated by the union of 
 the set ${\rm Unip} (G)$ of all unipotent elements in $G$, 
$D_{(\sigma_{i})_{s}}$ ($ 1\leqq i \leqq n$) and $H$. Then $M\cdot G^{0} = G.$
Consider a faithful finite-dimensional representation $\rho : G \to GL(V)$ and 
decompose $V = \oplus_{j = 1} ^{m} V_{\chi_{j}}$ for some distinct $\chi_{j} \in {\mathfrak X}(G^{0})$
with $V_{\chi_{j}} \not= \{0\}$.  As $G = Z_{G}(G^{0})$, the subspace $V_{\chi_{i}}$ is invariant under the action of $G$. The determinant ${\rm det}_{V_{\chi_{j}}}$  on $V_{\chi_{j}}$
induces a homomorphism 
$$\prod _{j=1}^{m}{\rm det}_{V_{\chi_{j}}} : G\ni \sigma  \mapsto ({\rm det}_{V_{\chi_{j}}} (\sigma))\in (\mbox{\boldmath $G_{m}$})^{m}$$
of $G$ to the product of $m$-copies of the multiplicative group $\mbox{\boldmath $G_{m}$}$
over $K$,  whose kernel is finite. Since $\prod _{j=1}^{m}{\rm det}_{V_{\chi_{j}}}({\rm Unip}(G))$
is trivial, $\prod _{j=1}^{m}{\rm det}_{V_{\chi_{j}}}(M)$ is a finitely generated
torsion subgroup of  $(\mbox{\boldmath $G_{m}$})^{m}$. Thus $M$ is the finite group
desired in this lemma.   \qed
    \begin{lem}\label{kernel}  Let $G$ be a connected reductive algebraic group such
   that $G = Z_{G}(T)$ for a closed connected subgroup $T$ of $G$
   which is an algebraic torus. Let $(X, G)$ be an effective regular action
   on an integral affine $K$-scheme $X = {\rm Spec} (R)$. 
   Suppose that ${\mathcal Q}(R)^{T} = {\mathcal Q} (R^{T})$. Then
   ${\rm Ker} (G \to {\rm Aut} (R^{T}))^{0} = T.$
   \end{lem}
   \proof 
   Put $H = {\rm Ker} (G \to {\rm Aut} (R^{T}))^{0}$.  Since the unipotent radical of $H$
   is normal in $G$, we see that $H$ is reductive.  Let $T'$ be a maximal
   connected torus of $H$.  As $T'\subseteq Z_{G}(T)$,  the subgroup $\widetilde{T} := T\cdot T'$ is  closed connected
    and diagonalizable in $G$. Hence 
   it  is an algebraic torus and $\widetilde{T} = T'$.  Consider the  group homomorphism
   $\rho : {\mathfrak X}(\widetilde{T}) \ni \chi \to \chi\vert_{T} \in {\mathfrak X}(T),$
   which is surjective. 
   For any $\chi_{i}\in {\mathfrak X}(\widetilde{T})^{R}$ let $w_{i}$ be a nonzero element 
   of $R_{\chi_{i}}$. 
   Then  for  almost all zero $a_{i} \in \mbox{\boldmath $Z$}$, 
   $$\prod_{i} w_{i}^{a_{i}} \in {\mathcal Q}(R)^{T} \Longleftrightarrow \sum_{i} a_{i}(\chi_{i}\vert_{T}) = 0.$$
   Since ${\rm Ker}(G \to {\rm Aut} (R^{T}) = {\rm Ker} (G \to {\rm Aut} ({\mathcal Q}(R)^{T})$
   contains $\widetilde{T}$, 
   from this equivalence we easily have 
   \begin{equation}\label{eqn2.5.1} \sum_{i} a_{i}\chi_{i} = 0 
    \Longleftrightarrow \sum_{i} a_{i}(\chi_{i}\vert_{T}) = 0. 
    \end{equation}
   If $\psi $ is a character in ${\mathfrak X}(T)^{R}$, decomposing a $\widetilde{T}$-module
    $R_\psi$ to  a direct sum
   of $R_\chi$  ($\chi\in {\mathfrak X}(\widetilde{T})^{R}$),  we see that 
   $\rho\vert _{\langle{\mathfrak X}(\widetilde{T})^{R}\rangle} : \langle{\mathfrak X}(\widetilde{T})^{R}\rangle \to \langle {\mathfrak X}({T})^{R}\rangle$
   is surjective. Thus by (\ref{eqn2.5.1}), $\rho\vert _{\langle{\mathfrak X}(\widetilde{T})^{R}\rangle}$
   is an isomorphism. As the action $(R, G)$ is effective, 
   both $\langle{\mathfrak X}(\widetilde{T})^{R}\rangle$ and 
   $\langle{\mathfrak X}({T})^{R}\rangle$ are of finite indices in ${\mathfrak X}(\widetilde{T})$
   and ${\mathfrak X}(T)$, respectively.  Consequently the kernel of $\rho$ is finite, which
   implies that ${\rm rank} (T) = {\rm rank} (\widetilde{T})$ and  $T= \widetilde{T}$. Since  any maximal torus of  the reductive $H$ is equal to $T$, we
   have  $H = T$ (e.g., \cite{Po}).  \qed 
   \begin{lem}\label{invarianttoricfield}
Let $(X, T)$ be a regular action of a connected
algebraic torus $T$ on  an integral affine $K$-scheme $X = {\rm Spec} (R)$. 
Then the following conditions are equivalent:
\begin{itemize}
\item[(i)]  ${\mathcal Q}(R)^{T} = {\mathcal Q} (R^{T})$.
\item[(ii)] $\dim _{{\mathcal Q}(R^{T})} {\mathcal Q}(R^{T})\otimes _{R^{T}} R_{\chi} =1$ for any
$ \chi \in {\mathfrak X}(T)^{R}$. 
\end{itemize}
\end{lem}
\proof
{\it (i)} $\Rightarrow$ {\it (ii)}:  Suppose that $R_{\chi}\not=\{0\}$ for a character $\chi\in {\mathfrak X}(T)$.  Let  $a, b \in R_{\chi}$ be any nonzero elements. As $\displaystyle\frac{a}{b}
\in {\mathcal Q}(R)^{T}$, by {\it (i)} we see
$a \in {\mathcal Q}(R^{T}) \cdot b$, which implies {\it (ii)}.

{\it (ii)} $\Rightarrow$ {\it (i)}: Let $x$ be any nonzero element of ${\mathcal Q}(R)^{T}$
and put  $ I_{x} := \left\{ b \in R \mid bx \in R\right\}.$  Then $I_{x}$
is a nonzero $T$-invariant ideal of $R$, which is a direct sum decomposition
$I_{x} = \bigoplus_{\psi \in {\mathfrak X}(T)} I_{x}\cap R_{\psi}$.   There is a
character $\chi \in {\mathfrak X}(T)$ such that $ I_{x}\cap R_{\chi}$ has a nonzero element $c$.
Since both $c$ and $cx$ is contained in $R_{\chi}$, by {\it (ii)} we can choose nonzero
$\alpha, \beta \in R^{T}$ in such a way that $\displaystyle\frac{\alpha}{\beta} c = cx$, 
which shows $x \in {\mathcal Q}(R^{T})$. 
 \qed


  \section{Toric Quotients and  Proof of  the Main Theorem} 
 
 Hereafter to the end of this paper, we suppose that $G$ is an algebraic group. 
 For  a regular action $(X, G)$
 with   $X= {\rm Spec} (R)$ and ${\mathfrak P}\in {\rm Ht}^{1}(R, R^{G^{0}})$, let 
 $\Delta_{G}({\mathfrak P})$ denote
 the $p$-part
 of  the order  of the factor group $$
\langle {\mathfrak X}(G^0)^{R^{{\mathcal I}_{G}({\mathfrak P})}}\rangle /\langle{\mathfrak X}(G^0)^{R^{{\mathcal I}_{G}({\mathfrak P})}/{\mathfrak P}\cap R^{{\mathcal I}_{G}({\mathfrak P})}}\rangle$$ if $p >0$. (Here the $p$-part $m_{p}$ of a natural number $m$ stands for 
the power of $p$ dividing $m$ such that $\displaystyle\frac{m}{m_{p}}$ is not divisible by
$p$.)  This order is finite, in  the  case of Theorem \ref{keytheorem1}.  If $p=0$, put  $\Delta_{G}({\mathfrak P})$=1. 
We say $(X, G)$ is {\it stable}, if  $X$ contains a non-empty open subset consisting of closed $G$-orbits
(cf. \cite{Po}).
 \begin{theorem}\label{keytheorem1}{\rm (cf. \cite{Nak2, Nak6}).}
Suppose that $G^0$ is an algebraic torus.  Then the following conditions are equivalent:
\begin{itemize}
\item[(i)]  $G = Z_G(G^0)$
\item[(ii)] For an arbitrary   closed subgroup $H$ of   $G$ containing $Z_G(G^0)$, 
the following 
conditions hold  for any 
effective regular action   $(X, H)$
on an arbitrary affine Krull $K$-scheme $X = {\rm Spec} (R)$ such that ${\mathcal Q}(R^{G^{0}}) = {\mathcal Q}(R)^{G^{0}}$:   for  any ${\mathfrak P}\in
{\rm Ht}_1(R, R^H)$, 
$${\rm e}({\mathfrak P}, {\mathfrak P}\cap R^H)={\rm e}({\mathfrak P}, {\mathfrak P}\cap R^{{\mathcal I}_H({\mathfrak P})}) \\
\cdot \Delta_{H}({\mathfrak P}).$$

\item[(iii)] For an arbitrary   closed subgroup $H$ of   $G$ containing $Z_G(G^0)$, 
the 
conditions in {\it (ii)} hold  for any 
effective stable  regular action   $(X, H)$
on an arbitrary affine normal variety   $X = {\rm Spec} (R)$.

\end{itemize}
\end{theorem}



\begin{pr}\label{unramified} Suppose that
$G^{0}$ is an algebraic torus. Let $(X, G)$ be an effective regular action of $G$
on  an affine Krull $K$-scheme $X = {\rm Spec} (R)$ such that  ${\mathcal Q}(R^{G^{0}}) =
{\mathcal Q} (R)^{G^{0}}$. Then we have 
$\langle {\mathcal I}_{G}({\mathfrak P})  \cup G^{0} \rangle =  {\mathcal I}_{G}({\mathfrak P}) \cdot G^{0} \subseteq Z_{G} (G^{0})$ and 
$${\rm e}({\mathfrak P}\cap R^{{\mathcal I}_{G}({\mathfrak P})\cdot G^{0}},
{\mathfrak  P}\cap R^{Z_{G}(G^{0})}) = 1$$
for any ${\mathfrak P} \in {\rm Ht}^{1}(R, R^{G}) = {\rm Ht}^{1}(R, R^{G^{0}})$.  
\end{pr}

\proof
Let ${\mathfrak P}$ be a prime ideal in  $ {\rm Ht}^{1}(R, R^{G})$.  Since  $R$
is a Krull domain, the orbit $G^{0} {\mathfrak P}$ is a finite set and hence  $\sigma ({\mathfrak P}) 
= {\mathfrak P}$ for any $\sigma \in G^{0}$ (e.g., \cite{Magid}). Then 
$$\sigma {\mathcal I}_{G}({\mathfrak P})\sigma^{-1}= 
 {\mathcal I}_{G}(\sigma({\mathfrak P})) =  {\mathcal I}_{G}({\mathfrak P}).$$
 As ${\mathcal I}_{G}({\mathfrak P})$ is a finite group (cf. \cite{Nak5}), we see ${\mathcal I}_{G}({\mathfrak P}) \subseteq Z_{G}(G^{0})$, which shows the first assertion.

 Put $H = {\mathcal I}_{G}({\mathfrak P}) \cdot G^{0}$. Then we have
   $H^{0} = G^{0}$, $H= Z_{H}(H^{0})$, 
   ${\mathcal I}_{H}({\mathfrak P}) = {\mathcal I}_{G}({\mathfrak P})$ and $\Delta_{G}({\mathfrak P})
   = \Delta_{H}({\mathfrak P})$. 
 Applying {\it (ii)} of Theorem \ref{keytheorem1} to   the induced action of $H$ on $X$, we
 must have
\begin{equation*} 
\begin{split} 
{\rm e}({\mathfrak P}, {\mathfrak P}\cap R^H)={\rm e}({\mathfrak P}, {\mathfrak P}\cap R^{{\mathcal I}_H({\mathfrak P})}) 
\cdot \Delta_{H}({\mathfrak P})\\
= {\rm e}({\mathfrak P}, {\mathfrak P}\cap R^{{\mathcal I}_G({\mathfrak P})}) 
 \cdot \Delta_{G}({\mathfrak P})
.\end{split}
\end{equation*} 
Replacing $H$ with  $H' = Z_{G}(G^{0})$,  by {\it (ii)} of Theorem \ref{keytheorem1} we similarly obtain
\begin{equation*} 
\begin{split} 
{\rm e}({\mathfrak P}, {\mathfrak P}\cap R^{H'})={\rm e}({\mathfrak P}, {\mathfrak P}\cap R^{{\mathcal I}_{H'}({\mathfrak P})}) 
 \cdot \Delta_{H'}({\mathfrak P})
\\
= {\rm e}({\mathfrak P}, {\mathfrak P}\cap R^{{\mathcal I}_G({\mathfrak P})}) 
 \cdot 
\Delta_{G}({\mathfrak P}).\end{split}
\end{equation*} 
Thus we have
${\rm e} ({\mathfrak P}\cap R^{H}, {\mathfrak P}\cap R^{H'} ) =1$.  \qed

\begin{pr}\label{keytheorem2}
Let $(X, T)$ be a regular action of a connected algebraic torus $T$
on an affine Krull $K$-scheme $X = {\rm Spec} (R)$.  Suppose that 
${\mathcal Q}(R)^{T} = {\mathcal Q}(R^{T})$.  Then $T$ acts naturally 
on $R/{\mathfrak P}$ satisfying 
$${\mathcal Q}(R/{\mathfrak P})^{T} =  {\mathcal Q}((R/{\mathfrak P})^{T}) $$
 for 
any ${\mathfrak P}\in {\rm Ht}^{1}(R, R^{T})$.  
\end{pr}
\proof
Fix a prime ideal ${\mathfrak P}\in {\rm Ht}^{1}(R, R^{T})$.  By the same reason as
in the proof of Proposition \ref{unramified}, the ideal  ${\mathfrak P}$
is invariant under the action of $T$ and the induced action $(R/{\mathfrak P}, T)$ 
is regular.  Put ${\mathfrak p} := {\mathfrak P}\cap R^{T}$, $U := R^{T}\backslash {\mathfrak P}$
and 
$A := R^{T}_{{\mathfrak p}} \otimes_{R^{T}} R$. We have naturally  a regular action 
$(A, T)$ of $T$ on a Krull $K$-domain $A$.  Note that
$(R/{\mathfrak P})^{T} = R^{T}/{\mathfrak p}$ and 
$${\mathcal Q} ((R/{\mathfrak P})^{T}) = {\mathcal Q}(R^{T}/{\mathfrak p}) = R^{T}_{\mathfrak p}/{\mathfrak p}R^{T}_{{\mathfrak p}}.$$ Let $\overline{x}$ (resp.  $\overline{U}$) denote
the image of $x\in R$ (resp. $U$) under the canonical morphism $R \to R/{\mathfrak P}$.  Suppose that 
$(R/{\mathfrak P})_{\chi}
\not= \{\overline{0}\}$ for a rational character $\chi \in {\mathfrak X}(T)$. Then
\begin{eqnarray*}
{\mathcal Q}((R/{\mathfrak P})^{T})\otimes _{(R/{\mathfrak P})^{T}} (R/{\mathfrak P})_{\chi}
\cong \overline{U}\/^{-1}(R/{\mathfrak P})_{\chi}. 
\end{eqnarray*}
On the other hand from a $T$-equivariant exact sequence
$$0 \to U^{-1}{\mathfrak P} \to A \to  \overline{U}\/^{-1}(R/{\mathfrak P}) \to 0$$ we have an exact sequence
\begin{equation*}
0 \to (U^{-1}{\mathfrak P})_\chi \to A_{\chi} \to  (\overline{U}\/^{-1}(R/{\mathfrak P}))_{\chi} \to 0
\end{equation*}
with canonical morphisms of $A^{T}$-modules. Thus $A_{\chi} \not= \{0\}$. 
Obviously $A^{T} \cong R^{T}_{{\mathfrak p}}$
and 
\begin{eqnarray}
 {\mathcal Q}(A^{T})\otimes_{A^{T}} A_{\chi} &\cong& {\mathcal Q}(R^{T}_{\mathfrak p})
\otimes _{R^{T}_{{\mathfrak p}}} R^{T}_{{\mathfrak p}} \otimes _{R^{T}} R_{\chi} \nonumber \\
&\cong& {\mathcal Q} (R^{T}) \otimes_{R^{T}} R_{\chi}.\nonumber\end{eqnarray}
As ${\mathcal Q}(A)^{T}\cong {\mathcal Q}(R)^{T} = {\mathcal Q} (R^{T})
\cong {\mathcal Q}(A^{T})$, by Lemma \ref{invarianttoricfield}  we see that
$$\dim_{{\mathcal Q}(A^{T})} {\mathcal Q}(A^{T}) \otimes_{A^{T}} A_{\chi} = 1.$$
Since $(R/{\mathfrak P})_{\chi} \not={\overline{0}}$, we can choose an element $f$ from $R_{\chi}$
in such a way that $\overline{f} \not=\overline{0}$. 
Clearly 
${\rm Over}_{{\mathfrak p}A^{T}}(A) = \{ {\mathfrak Q}A \mid {\mathfrak Q}\in  {\rm Over}_{\mathfrak p}(R)\}$
and $$0 = {\rm v}_{A, {\mathfrak P}A}(1\otimes f) < {\rm e}({\mathfrak P}A, {\mathfrak p}A^{T})
= {\rm e}({\mathfrak P}, {\mathfrak p}).$$
Thus by Theorem 3.3 of \cite{Nak3}, we must have
$A_{\chi} = A^{T}\cdot 1\otimes f,$
which implies that $$\dim _{{\mathcal Q}((R/{\mathfrak P})^{T})} {\mathcal Q}((R/{\mathfrak P})^{T})\otimes _{(R/{\mathfrak P})^{T}} (R/{\mathfrak P})_{\chi} = 1.$$
Consequently from Lemma \ref{invarianttoricfield} we establish the assertion of this proposition. \qed

\begin{theorem}\label{main}  Let   $G^{0}$ be an algebraic torus and     
$(X, G)$  an effective regular action of $G$ on an affine Krull $K$-scheme
$X = {\rm Spec} (R)$ such that ${\mathcal Q}(R)^{G^{0}}= {\mathcal Q}(R^{G^{0}})$. 
Let $N$ be any closed normal subgroup of $G$ containing $G^{0}$.  If  $G = Z_{G}(G^{0})$, then 
$$({\mathcal I}_{G}({\mathfrak P})\vert_{R^{N}} = ({\mathcal I}_{G}({\mathfrak P})\cdot N)\vert _{R^{N}}
= {\mathcal I}_{G}({\mathfrak P}\cap R^{N})\vert _{R^{N}}$$
for any prime ideal ${\mathfrak P}\in  {\rm Ht}^{1}(R, {G^{0}})$ and the following equality
holds: 
$${\mathfrak R} (R^{N}, G)\vert _{R^{N}}= {\mathfrak R} (R, G) \vert_{R^{N}}.$$
 
\end{theorem}

\proof Suppose that   $N = G^{0}$. 
 Let ${\mathfrak P}$ be any prime ideal in ${\rm Ht}^{1}(R, G^{0})
= {\rm Ht}^{1}(R, G)$ and put ${\mathfrak p} := {\mathfrak P}\cap R^{G^{0}}$. 
As  $${\mathcal I}_{G}({\mathfrak p}) = \{ \sigma \in G \mid (\sigma - 1)(R^{G^{0}})
\subseteq {\mathfrak p} \},$$ we see ${\mathcal I}_{G}({\mathfrak P})\cdot G^{0} \subseteq {\mathcal I}_{G}({\mathfrak p})$. 
By Lemma \ref{liftpart} we can choose  a finite subgroup $M$ of $G$ in such a way that 
$M \supseteq {\mathcal I}_{G}({\mathfrak P})$ and 
\begin{equation}\label{eqn:theorem}  M \cdot G^{0} = {\mathcal I}_{G}({\mathfrak p}).\end{equation}  Then from definition of $ {\mathcal I}_{G}({\mathfrak P})$ we must have 
${\mathcal I}_{M}({\mathfrak P})
= {\mathcal I}_{G}({\mathfrak P}).$
Consequently by Chap. VI of \cite{LR}  the field 
${\mathcal Q}(R^{{\mathcal I}_{G}({\mathfrak P})}/{\mathfrak P}\cap R^{{\mathcal I}_{G}({\mathfrak P})})$
is a Galois extension over the field
$${\mathcal Q}(R^{{\mathcal D}_{M}({\mathfrak P})}/{\mathfrak P}\cap R^{{\mathcal D}_{M}({\mathfrak P})}) = {\mathcal Q}(R^M/{\mathfrak P}\cap R^M)$$ 
with the Galois group ${\mathcal D}_{M}({\mathfrak P})/{\mathcal I}_{M}({\mathfrak P})$
(for ${\mathcal D}_{M}({\mathfrak P})$, cf.  {\bf 1A}). 
From Proposition \ref{keytheorem2} we infer that
\begin{eqnarray*} 
{\mathcal Q}(R^{{\mathcal I}_{G}({\mathfrak P})\cdot G^{0}}/{\mathfrak P}\cap R^{{\mathcal I}_{G}({\mathfrak P})\cdot G^{0}})  = {\mathcal Q}(R^{{\mathcal I}_{G}({\mathfrak P})}/{\mathfrak P}\cap R^{{\mathcal I}_{G}({\mathfrak P})})^{G^{0}}\end{eqnarray*}
and 
\begin{eqnarray*}
  {\mathcal Q}(R^{{\mathcal D}_{M}({\mathfrak P})\cdot G^{0}}/{\mathfrak P}\cap R^{{\mathcal D}_{M}({\mathfrak P})\cdot G^{0}}) 
& = & {\mathcal Q}(R^{{\mathcal D}_{M}({\mathfrak P})}/{\mathfrak P}\cap R^{{\mathcal D}_{M}({\mathfrak P})})^{G^{0}} \\ 
&=&  {\mathcal Q}(R^{M}/{\mathfrak P}\cap R^{M})^{G^{0}} \\
& = &
{\mathcal Q}(R^{M\cdot G^{0}}/{\mathfrak P}\cap R^{M\cdot G^{0}}). 
\end{eqnarray*}
On the other hand 
\begin{eqnarray*} {\mathcal Q}(R^{{\mathcal D}_{M}({\mathfrak P})}/{\mathfrak P}\cap R^{{\mathcal D}_{M}({\mathfrak P})})^{G^{0}} 
 &=& ({\mathcal Q}(R^{{\mathcal I}_{G}({\mathfrak P})}/{\mathfrak P}\cap R^{{\mathcal I}_{G}({\mathfrak P})})^{{\mathcal D}_{M}({\mathfrak P})/{\mathcal I}_{M}({\mathfrak P})})^{G^{0}} \\
&= & ({\mathcal Q}(R^{{\mathcal I}_{G}({\mathfrak P})}/{\mathfrak P}\cap R^{{\mathcal I}_{G}({\mathfrak P})})^{G^{0}})^{{\mathcal D}_{M}({\mathfrak P})/{\mathcal I}_{M}({\mathfrak P})}.\end{eqnarray*}
Consequently  the field ${\mathcal Q}(R^{{\mathcal I}_{G}({\mathfrak P})\cdot G^{0}}/{\mathfrak P}\cap R^{{\mathcal I}_{G}({\mathfrak P})\cdot G^{0}})$ is a finite Galois extension
of the field $${\mathcal Q}(R^{M\cdot G^{0}}/{\mathfrak P}\cap R^{M\cdot G^{0}})
= {\mathcal Q}(R^{{\mathcal I}_{G}({\mathfrak p})} /{\mathfrak p}\cap 
R^{{\mathcal I}_{G}({\mathfrak p})})$$ 
(cf. (\ref{eqn:theorem})).
However by Chap. VI of \cite{LR}  the field of ${\mathcal Q} (R^{G^{0}}/{\mathfrak p})$ is purely inseparable over the field 
$${\mathcal Q}(R^{{\mathcal I}_{G}({\mathfrak p})} /{\mathfrak p}\cap 
R^{{\mathcal I}_{G}({\mathfrak p})}).$$ 
  Hence we must have
$${\mathcal Q}(R^{{\mathcal I}_{G}({\mathfrak P})\cdot G^{0}}/{\mathfrak P}\cap R^{{\mathcal I}_{G}({\mathfrak P})\cdot G^{0}})= {\mathcal Q}(R^{{\mathcal I}_{G}({\mathfrak p})} /{\mathfrak p}\cap 
R^{{\mathcal I}_{G}({\mathfrak p})})$$
which shows the discrete valuation rings  $$A:=(R^{{\mathcal I}_{G}({\mathfrak P})\cdot G^{0}})_{{\mathfrak P}\cap R^{{\mathcal I}_{G}({\mathfrak P})\cdot G^{0}}}, 
B:= (R^{{\mathcal I}_{G}({\mathfrak p})})_{{\mathfrak p}\cap 
R^{{\mathcal I}_{G}({\mathfrak p})}}$$  
with maximal ideals ${\mathfrak M}_{A}$, ${\mathfrak M}_{B}$ respectively
have the common residue class field, i.e., $A/{\mathfrak M}_{A} = B/{\mathfrak M}_{B}$. The field 
${\mathcal Q}(R^{G^{0}})$ is a  Galois extension over ${\mathcal Q}(B) =
{\mathcal Q}(R^{G^{0}})^{{\mathcal I}_{G}({\mathfrak p})/G^{0}}$ and the ideal 
${\mathfrak p}$ is invariant under the action of ${\mathcal I}_{G}({\mathfrak p})$. 
Hence   $A$ 
 is the  integral closure
of  $B$ 
 in a finite separable field extension, which implies 
$A$ 
is a finite 
$B$-module (e.g., \cite{BourbakiC}).  On the other hand we infer from 
Proposition \ref{unramified} that  ${\mathfrak M}_{A} = A\cdot {\mathfrak M}_{B}$. 
Consequently we have   $A = B$, which shows  
${\mathcal Q}(A) = {\mathcal Q}(B)$.   By Galois theory
we see $${\mathcal I}_{G}({\mathfrak P}) \vert_{R^{G^{0}}}= ({\mathcal I}_{G}({\mathfrak P})\cdot G^{0})\vert _{R^{G^{0}}}
= {\mathcal I}_{G}({\mathfrak p})\vert _{R^{G^{0}}}.$$
The last assertion ${\mathfrak R} (R^{G^{0}}, G)\vert _{R^{G^{0}}}= {\mathfrak R} (R, G) \vert_{R^{G^{0}}}$ of the theorem 
follows immediately from this equality 
and ${\rm Over}_{\mathfrak q}(R)\not=\emptyset$ for any  ${\mathfrak q}\in{\rm Spec} (R^{G^{0}})$
of ${\rm ht}({\mathfrak q}) = 1$ (e.g. \cite{Magid}). 

Now we treat the case of $N\varsupsetneq G^{0}$. 
Appling the proof of Proposition \ref{finitecase}  to the action $(R^{G^{0}}, G/G^{0})$, 
we see that 
$$({\mathcal I}_{G}({\mathfrak p})\cdot N)\vert _{R^{N}}
= ({\mathcal I}_{G/G^{0}}({\mathfrak p})\cdot N/G^{0})\vert _{R^{N}}
= {\mathcal I}_{G}({\mathfrak p}\cap R^{N})\vert _{R^{N}}.$$
This and the first paragraph of the proof imply 
\begin{eqnarray*} ({\mathcal I}_{G}({\mathfrak P})\cdot N)\vert _{R^{N}}
&=& ({\mathcal I}_{G}({\mathfrak P})\vert_{R^{G^{0}}}\cdot N\vert_{R^{G^{0}}})\vert _{R^{N}}\\
&= &({\mathcal I}_{G}({\mathfrak p})\vert _{R^{G^{0}}}\cdot N\vert_{R^{G^{0}}})\vert _{R^{N}} 
= ({\mathcal I}_{G}({\mathfrak p})\cdot N)\vert _{R^{N}} \\
&=&{\mathcal I}_{G}({\mathfrak P}\cap R^{N})\vert _{R^{N}}\end{eqnarray*}
and the last assertion of the theorem follows from the equality. 
\qed

   \noindent {\bf Proof of Theorem \ref{maintheorem}~:}  Clearly 
   $${\mathfrak R}(R, G ; T)\cdot T\subseteq {\mathfrak R}(R^{T}, G).$$ Since 
   ${\mathfrak R}(R^{T}, G)\vert_{R^{T}}$ is finite, by Lemma \ref{kernel}, the torus 
   $T$ is  the identity
   component of the closed subgroup $H : = {\mathfrak R}(R^{T}, G)$ of $G$. 
   Then ${\mathfrak R}(R, G;  T) = {\mathfrak R}(R, H)$ and ${\mathfrak R}(R^{T}, G)
   = {\mathfrak R}(R^{T}, H)$. 
  Thus the  equality  $${\mathfrak R}(R, G;  T)\vert_{R^{T}}= {\mathfrak R}(R^{T}, G)\vert_{R^{T}}$$ 
  follows from Theorem \ref{keytheorem2} for $N= H^{0}= T$. 
  \qed

    
    \small \section{Comments}\normalsize
    
   
  {\bf 4.A.}   Concerning Theorem \ref{main}, we show that the assumption $G = Z_{G}(G^{0})$
    is indispensable for the assertion as follows: 
    \begin{pr}\label{noncentral} Suppose that
$G^{0}$ is an algebraic torus and $G \not= Z_{G}(G^{0})$. 
Then there exists a closed subgroup $H$ of $G$ with
$H \supsetneqq Z_{G}(G^{0})$ and an effective  regular action $(X, H)$ of
$H$ on an affine normal variety $X = {\rm Spec} (R)$ such that    ${\mathcal Q}(R^{G^{0}}) =
{\mathcal Q} (R)^{G^{0}}$ and 
$${\rm e}({\mathfrak P}\cap R^{ {\mathcal I}_{H}({\mathfrak P})\cdot G^{0}},
{\mathfrak  P}\cap R^{{\mathcal I}_{H}({\mathfrak P}\cap R^{G^{0}})}) > 1$$
for some ${\mathfrak P} \in {\rm Ht}_{1}(R, R^{G})$.  Moreover we have
$${\mathfrak R}(R, H) \vert_{R^{G^{0}}} \not= {\mathfrak R}(R^{G^{0}}, H)\vert_{R^{G^{0}}}.$$
\end{pr}

\proof By {\it (iii)} of Theorem \ref{keytheorem1}, there is a  closed subgroup $H$ of $G$ with
$H \supsetneqq Z_{G}(G^{0})$ and an effective  regular action $(X, H)$ of
$H$ on an affine normal variety $X = {\rm Spec} (R)$ such that    ${\mathcal Q}(R^{G^{0}}) =
{\mathcal Q} (R)^{G^{0}}$ and 
\begin{equation}\label{nonequality}{\rm e}({\mathfrak P}, {\mathfrak P}\cap R^H)\not={\rm e}({\mathfrak P}, {\mathfrak P}\cap R^{{\mathcal I}_H({\mathfrak P})}) 
\cdot \Delta_{H}({\mathfrak P}).\end{equation}
(The stability of $(X, H)$ implies ${\mathcal Q}(R^{G^{0}}) =
{\mathcal Q} (R)^{G^{0}}$.) 
From the proof of Proposition \ref{unramified} we infer that 
the right hand side of (\ref{nonequality}) is equal to ${\rm e}({\mathfrak P}, {\mathfrak P}\cap R^{Z_{H}(G^{0})})$, which implies 
\begin{equation}\label{eqn:noncentral} 
{\rm e}({\mathfrak P}, {\mathfrak P}\cap R^H)  > {\rm e}({\mathfrak P}, {\mathfrak P}\cap R^{Z_{H}(G^{0})})  =  {\rm e}({\mathfrak P}, {\mathfrak P}\cap R^{{\mathcal I}_{H}({\mathfrak P})\cdot G^{0}}).\end{equation}
Since $H \supseteq {\mathcal I}_{H}({\mathfrak P}\cap R^{G^{0}})\supseteq G^{0}$, by {\it (ii)} of Lemma \ref{normalcase} we see $
{\rm e}({\mathfrak P}\cap R^{ {\mathcal I}_{H}({\mathfrak P}\cap R^{G^{0}})}, {\mathfrak P}\cap R^H) =1.$
 Consequently the inequality of ramification indices of  Proposition \ref{noncentral} follows from 
 (\ref{eqn:noncentral}).  By this we see that 
 \begin{equation}\label{eqn:nonequ} ({\mathcal I}_{H}({\mathfrak P})\cdot G^{0})\vert _{R^{G^{0}}}\subsetneq
 {\mathcal I}_{H}({\mathfrak P}\cap R^{G^{0}})\vert _{R^{G^{0}}}.
 \end{equation}
 As ${\mathcal I}_{H}({\mathfrak P}\cap R^{G^{0}}) \subseteq {\mathfrak R}(R^{G^{0}}, H)$, 
 we have \begin{equation}\label{eqn:4.1.1} {\mathcal I}_{H}({\mathfrak P}\cap R^{G^{0}}) =
 {\mathcal I}_{{\mathfrak R}(R^{G^{0}}, H)}({\mathfrak P}\cap R^{G^{0}}).
 \end{equation}
 Similarly we see ${\mathcal I}_{H}({\mathfrak P}) =
 {\mathcal I}_{{\mathfrak R}(R, H)}({\mathfrak P})$ and hence
 \begin{equation}\label{eqn:4.1.2}
 {\mathcal I}_{H}({\mathfrak P}) =  {\mathcal I}_{{\mathfrak R}(R, H)\cdot G^{0}}({\mathfrak P}).
 \end{equation}
 
 Assume 
 $${\mathfrak R}(R, H)\cdot G^{0} \vert_{R^{G^{0}}} = {\mathfrak R}(R^{G^{0}}, H)\vert_{R^{G^{0}}}.$$
 Since ${\mathfrak R}(R, H)$ is a finite group contained in $Z_{H}(G^{0})$, applying Theorem \ref{main}
 to the action $(X, {\mathfrak R}(R, H)\cdot G^{0})$, we see 
 $${\mathcal I}_{{\mathfrak R}(R, H)\cdot G^{0}}({\mathfrak P})\cdot G^{0}
\vert _{R^{G^{0}}}={\mathcal I}_{{\mathfrak R}(R, H)\cdot G^{0}}({\mathfrak P}\cap R^{G^{0}})
\vert _{R^{G^{0}}} = {\mathcal I}_{{\mathfrak R}(R^{G^{0}}, H)}({\mathfrak P}\cap R^{G^{0}})
\vert _{R^{G^{0}}} .$$
This equality, (\ref{eqn:4.1.1}) and  (\ref{eqn:4.1.2}) imply $${\mathcal I}_{H}({\mathfrak P})\cdot G^{0}\vert _{R^{G^{0}}}=
 {\mathcal I}_{H}({\mathfrak P}\cap R^{G^{0}})\vert _{R^{G^{0}}},$$
 which conflicts with (\ref{eqn:nonequ}).  Hence  we must have the last assertion
 of this proposition. \qed 
 
 \medskip
\noindent{\bf 4.B.} 
Concerning Proposition \ref{keytheorem2}  which is fundamental in the proof of our theorem, 
we give the following remark: 
\begin{rem}  
 The condition ``${\mathcal Q}(R)^{T}= {\mathcal Q}(R^{T})$''  does not always imply 
 ``${\mathcal Q}(R/{\mathfrak P})^{T}= {\mathcal Q}((R/{\mathfrak P})^{T})$'' 
 in the case that ${\rm ht}({\mathfrak P}) = 1 < {\rm ht}({\mathfrak P}\cap R^{T})$. 
  We give an example as follows: 
\rm Let $V: = K^{5}$  and 
\begin{eqnarray*} T: =\left\{\left[
\begin{array}{ccccc} 
t &  & & &  \\  & t^{-1} & & & \\ 
& & t^{-1} &  & \\ & & & u  & \\
& & & &  u^{-1}\end{array} \right]  \mid  (t, u) \in (\mbox{\boldmath $G_{m}$})^{2} \right\} 
\end{eqnarray*}  where the matrix representation on the dual space $V^{\vee}$
of $V$ defined by the basis $$\{X_{1}, X_{2}, X_{3}, X_{4}, X_{5}\}.$$
Let $R :=K[X_{1}, X_{2}, X_{3}, X_{4}, X_{5}]$ be a $5$-dimensional polynomial 
ring on which $T$ acts naturally and put ${\mathfrak P}:= R X_{1}$. By stability
of $(R, T)$, we have 
${\mathcal Q}(R)^{T} = {\mathcal Q} (R^{T})$. Moreover
${\mathfrak P}\cap R^{T} = R^{T}( X_{1}X_{2},  X_{1}X_{3})$
and $(R/{\mathfrak P})^{T} = K[\overline{X_{4}}\cdot \overline{X_{5}}]$. 
Here  $\overline{x} := x +{\mathfrak P}\in R/{\mathfrak P}$.  Thus
${\rm ht} ({\mathfrak P}) = 1 < 2 ={\rm ht} ({\mathfrak P}\cap R^{T})$ and
$${\mathcal Q}(R/{\mathfrak P})^{T}\ni \frac{\overline{X_{2}}}{\overline{X_{3}}}
\not\in {\mathcal Q}((R/{\mathfrak P})^{T}).$$
\end{rem}
  
  \medskip
\noindent  {\bf 4.C.}  Concerning Theorem \ref{main}, we note the two remarks.  
   \begin{rem}\label{simple} The equality ${\mathfrak R}(G, R)\cdot G^{0}
   \vert _{R^{G^{0}}}
   = {\mathfrak R}(G, R^{G^{0}})\vert _{R^{G^{0}}}$ does not necessarily  hold  for
   a semi-simple $G^{0}$. We give an example to explain this as follows: \rm
   Let $G^{0}$ be ${ SL}_{n}$ with $n \geqq 2$ and let 
   $\Phi _{1}$ be the  $n$-dimensional representation on which $G^{0}$ is standard. 
   Let $V$ be an $n^{2}$-dimensional 
   representation
   $$V = \underbrace{\Phi_{1} \oplus \cdots \oplus \Phi_{1}}_{\text{$n$}}$$
   Let $\sigma \in GL(V)$ be the
    scalar matrix $\zeta_{d}\cdot E_{n^{2}}$, where $E_{n^{2}}$ is the unit matrix in $GL(V)$
    and $\zeta_{d}$ is a fixed primitive $d$-th root of $1\in K$ with
    a natural number $d> n$ such that $p$ does not divide $d$ if $p>0$.  Define
    $G$ to be the subgroup of $GL(V)$ generated by $G^{0}$ in $GL(V)$ and $\sigma$. 
    Let $R$ be a $K$-algebra of polynomial functions on $V$. Then $R^{G^{0}}  = K[g]$
    for a homogeneous polynomial $g$ of degree $n$ and hence 
    $${\mathfrak R} (R^{G^{0}}, G)\vert _{R^{G^{0}}} = \langle \sigma\rangle \vert_{R^{G^{0}}},$$
    which is a non-trivial group.  On the other hand, as ${\mathfrak R} (R, G)$ is finite
    and $\sigma$ is a scalar matrix,
    ${\mathfrak R} (R, G)$ is generated by pseudo-reflections which are 
    scalar matrices in $GL(V)$.  Thus
    $$\{1\} = {\mathfrak R}(R, G)\vert _{R^{G^{0}}} \not= 
    {\mathfrak R} (R^{G^{0}}, G)\vert _{R^{G^{0}}}.$$

   \end{rem}
   
      \begin{rem}\label{commoncounter} 
   The assumption ${\mathcal Q}(R^{G^{0}}) = {\mathcal Q}(R)^{G^{0}}$ is indispensable 
   for  Theorem \ref{main}. We show an example as follows: 
   \rm  Let $\zeta_{d}$ be the same element as in Remark \ref{simple}. 
    Let $V: = K^{4}$  and define $G = \langle T, \sigma \rangle\subseteq GL(V^{\vee})$ by 
\begin{eqnarray*} T: =\left\{\left[
\begin{array}{ccccc} 
t  &  & &   \\  & t^{-1} & &  \\ 
& & u  &   \\ & & & u   
\end{array} \right]  \mid  (t, u) \in (\mbox{\boldmath $G_{m}$})^{2} \right\} ,
\sigma := \left[
\begin{array}{ccccc} 
\zeta_{d }&  & &   \\  & 1 & &  \\ 
& & \zeta_{d} &   \\ & & & 1  
\end{array} \right] 
\end{eqnarray*}  where the matrix representation is given on the basis $\{X_{1}, X_{2}, X_{3}, X_{4}\}$ of the dual space $V^{\vee}$   of $V$.  Clearly $G^{0} = T$ and $G = Z_{G}(G^{0})$. Let $R$ be the  $K$-algebra of polynomial functions on
$V$ and ${\mathfrak P} = R X_{1}$.  Then ${\mathcal I}_{G}({\mathfrak P}) = \{1\}$.   Since $R^{G^{0}} = K[X_{1}X_{2}]$, 
we see ${\mathfrak P} \in {\rm Ht}_{1}(R, R^{G^{0}})$ and ${\mathcal I}_{G}({\mathfrak P}\cap R^{G^{0}}) = G$. Thus $$\{1\} = ({\mathcal I}_{G}({\mathfrak P}) \cdot G^{0})\vert_{R^{G^{0}}} \not=
 {\mathcal I}({\mathfrak P}\cap R^{G^{0}}) \vert_{R^{G^{0}}} \cong  G/G^{0}\cong \langle \sigma
 \rangle.$$  On the other hand 
 $${\mathcal Q}(R)^{G^{0}} = K\left(X_{1}X_{2}, \frac{X_{3}}{X_{4}}\right)\not= 
 {\mathcal Q}(R^{G^{0}}).$$
   \end{rem} 
   
   \medskip
 \noindent {\bf 4.D.} 
    For a rational $G^{0}$-module $M$, let
          ${\mathfrak X}(G^{0})_{M} = \{ \chi \in {\mathfrak X}(G^{0}) \mid 
        M_{\chi}\not=\{0\}, M_{-\chi} \not= \{0\} \}$ and put ${\mathcal N}(G^{0}, M)=  \bigcap_{\chi \in {\mathfrak X}(G^{0})_{M}} {\rm Ker} (\chi)$. 
        
        \begin{rem}\label{maximalstable}
       Suppose that $G$ is an algebraic connected torus and let $(X, G)$ be an effective 
       regular action 
       of $G$ on an  affine $K$-scheme $X = {\rm Spec} (R)$.  \rm 
        By the definition of ${\mathfrak   X}(G)_{R}$,
        we easily see that ${\mathfrak   X}(G/{\mathcal  N}(G, R)) \cong {\mathfrak X}(G)_{R}$. 
        Now suppose that {\it $R$ is finitely generated as a $K$-algebra.} Then
         $(R^{{\mathcal  N}(G, R)}, G)$ is the  stable action of $G$ on the largest $K$-subalgebra
         in  $R$. \it
Consider  the following conditions on $(X, G)$ :
       \begin{itemize}
       \item[(i)] $(X, G)$ is a stable action. 
       
       \item[(ii)] ${\mathfrak X}(G)^{R} = {\mathfrak X}(G)_{R}$.
       
      \item[(iii)] ${\mathcal Q}(R)^{G} = {\mathcal Q}(R^{G})$.
        
        \end{itemize}
        Then (i) $\Leftrightarrow$ (ii) $\Rightarrow$ (iii) are easy. 
        (ii) $\Rightarrow$ (iii)  holds in the case where $R$ is not necessarily finitely generated
        as a $K$-algebra.  However the implication (ii) $\Leftarrow$ (iii)
        is not true.  \rm  In Remark \ref{commoncounter} the 
         action $(K[X_{1}, X_{2}, X_{3}], G^{0})$
        is not stable, but ${\mathcal Q}(K[X_{1}, X_{2}, X_{3}]^{G^{0}})=
        K(X_{1}, X_{2}, X_{3})^{G^{0}}.$        
        \end{rem}

Since ${\mathcal Q}(R^{{\mathcal N}(T, R)})^{T} = {\mathcal Q}(R^{T})$
(cf.  Remark \ref{maximalstable})
and $R^{{\mathcal N}(T, R)}$ is a Krull $K$-domain with a regular action of $G$
in the following case, 
we immediately have a version of Theorem \ref{maintheorem}
without the assumption that ${\mathcal Q}(R)^{T} = {\mathcal Q}(R^{T})$:

\begin{theorem}
Suppose that $G^0$ is reductive and let $T$ be 
 a connected closed subgroup of $G^{0}$ which is   an algebraic torus. 
Let $(X, G)$ be an effective  regular action of $G$ on an affine Krull 
$K$-scheme $X$ with $X = {\rm Spec} (R)$. Suppose that $G = Z_{G}(T)$. 
 Then we have  
$${\mathfrak R}(R^{{\mathcal N}(T, R)}, G; T)\vert _{R^{T}} = {\mathfrak R}(R^{T}, G)\vert _{R^{T}}.
~\mbox{\qed}$$ 
\end{theorem}

\medskip
\noindent {\bf 4.E.}   In Problem 1.1 we may suppose that  the pseudo-reflection group ${\mathfrak R}(R, G)$
         is finite on $R$.  This finiteness condition holds
         for all  regular actions
           $(R, G)$ of $G$ on Krull  $K$-domains  if and only if $G^{0}$ is reductive (cf. \cite{Nak5}). 
           Thus the ``reductive'' assumption of $G^{0}$ is necessary for the main theorem in this case.

    
    \small \section{Chevalley-Serre and Steinberg  Theorems}\normalsize
    

For a finite dimensional faithful representation $H \to GL(V)$ of a finite
group $H$ over $K$, if $R^{H}$ is a polynomial ring over $K$, then $H$ is generated by pseudo-reflections  (i.e., $H = {\mathfrak R}(R, H)$)
where $R$ denotes the $K$-algebra of polynomial functions on $V$. 
The converse
of the assertion  is true, especially in the case where $p=0$ or the order of $H$ is not divisible by $p$.  These results  are obtained by G.C. Shephard, J.A. Todd, 
C. Chevalley and J.-P. Serre (e.g.,  \cite{Bourbaki, Serre}). 

 We say
a commutative ring $A$ is regular (resp.  a locally complete intersection), if  all 
 localizations of $A$ at prime ideals are regular
(resp.  complete intersections).  Let ${\rm m\mathchar`-Spec} (A)$ denote the
 maximum spectrum of a commutative ring $A$.  Clearly,  for 
a finitely generated positively graded algebra $A$ defined  over $K$, $A$ is a polynomial ring
over $K$ (resp. a complete intersection) if and only if $A_{A_{+}}$ is regular (resp. so), where $A_+$ is the 
homogeneous maximal ideal of $A$.  We 
     generalize the first  assertion 
     and a similar result obtained by V.G. Kac and  K.-I. Watanabe (\cite{Serre, Kac})
     on complete intersections 
as follows.  
\begin{pr}\label{generalSerre}  Let $S$ be an integrally closed domain and  $H$ a finite subgroup of 
${\rm Aut} (S)$. 
Suppose that ${\mathcal D}_{H}({\mathfrak M}) = {\mathcal I}_{H}({\mathfrak M})$  
 for any  ${\mathcal D}_{H}({\mathfrak M})$ which is  maximal in 
$$\{ {\mathcal D}_{H}(\tilde {\mathfrak M}) \mid \tilde {\mathfrak M}\in {\rm m \mathchar`- Spec} (S) \}$$
with respect to inclusions. Let $N$ be a normal subgroup of $H$
containing all  ${\mathcal I}_{H}({\mathfrak P})$  for  prime ideals  ${\mathfrak P}$
of $S$ of height $1$ (resp. of height $\leqq 2$).  If $S^{H}$ is regular (resp. a locally complete intersection), then 
${\mathcal D}_{H}({\mathfrak m}) \subseteq N$ for any maximal ideal ${\mathfrak m}$
of $S^{N}$. \end{pr}

\proof  We treat this in  the case where $S^{H}$ is regular. 
Let ${\mathfrak M}$ be any maximal ideal of $S$. First, we will show
\begin{equation}\label{inclusion}  {\mathcal D}_{H}({\mathfrak M}) \subseteq N. \end{equation}
To show this we may assume that   ${\mathcal D}_{H}({\mathfrak M})$  is a maximal   in 
$\{ {\mathcal D}_{H}(\tilde{\mathfrak M}) \mid \tilde{\mathfrak M}\in 
{\rm m\mathchar`-Spec} (S) \}.$ Since 
$(S^{{\mathcal D}_{H}({\mathfrak M})})_{{\mathfrak M}\cap S^{{\mathcal D}_{H}({\mathfrak M})} }$ 
is \'etale over $S^{H}_{{\mathfrak M}\cap S^{H}}$ (cf. (1) of Proposition \ref{finitecase}), the local 
ring 
$(S^{{\mathcal D}_{H}({\mathfrak M})})_{{\mathfrak M}\cap S^{{\mathcal D}_{H}({\mathfrak M})} }$
is regular.  Let ${\mathfrak p}$ be any prime ideal of $S^{{\mathcal D}_{N}({\mathfrak M})}$ of height $1$ contained in ${\mathfrak M}\cap S^{{\mathcal D}_{N}({\mathfrak M})}$ and 
${\mathfrak P}$ a prime ideal in ${\rm Over}_{{\mathfrak p}}(S)$ contained in
${\mathfrak M}$. 
From the definition of $N$, we see  
${\mathcal I} _{{\mathcal D}_{H}({\mathfrak M})}({\mathfrak P})
\subseteq {\mathcal D}_{N}({\mathfrak M})= N \cap  {\mathcal D}_{H}({\mathfrak M})$ and moreover,  by (1) of Proposition \ref{finitecase},  we 
see $(S^{{\mathcal D}_{N}({\mathfrak M})})_{\mathfrak p}$ is \'etale over
$(S^{{\mathcal D}_{H}({\mathfrak M})})_{{\mathfrak p}\cap S^{{{\mathcal D}_{H}({\mathfrak M})}}}$. 
As $(S^{{\mathcal D}_{N}({\mathfrak M})})_{{\mathfrak M}\cap S^{{\mathcal D}_{N}({\mathfrak M})} }$ 
is the integral closure of 
$(S^{{\mathcal D}_{H}({\mathfrak M})})_{{\mathfrak M}\cap S^{{\mathcal D}_{H}({\mathfrak M})} }$ 
in a finite separable (Galois) extension of quotient fields, 
$(S^{{\mathcal D}_{N}({\mathfrak M})})_{{\mathfrak M}\cap S^{{\mathcal D}_{N}({\mathfrak M})} }$
is noetherian.  Since 
the local homomorphism \begin{equation}\label{serre} (S^{{\mathcal D}_{H}({\mathfrak M})})_{{\mathfrak M}\cap S^{{\mathcal D}_{H}({\mathfrak M})} } \to (S^{{\mathcal D}_{N}({\mathfrak M})})_{{\mathfrak M}\cap S^{{\mathcal D}_{N}({\mathfrak M})} }\end{equation}
 is unramified of codimension $1$, applying Nagata's purity 
of branch loci (cf. \cite{LR}) to the inclusion, we infer that
this morphism is unramified. Thus the induced morphism  of (\ref{serre}) of residue class fields is regarded as  a Galois extension under the
action of the group  
${\mathcal D}_{H}({\mathfrak M})/{\mathcal D}_{N}({\mathfrak M})$ (e.g., Chap. V, \cite{LR}). 
By our assumption of this proposition, the local rings of (\ref{serre})
have the common residue class field and, as (\ref{serre}) is finite, 
by Nakayama's lemma we must have $$(S^{{\mathcal D}_{H}({\mathfrak M})})_{{\mathfrak M}\cap S^{{\mathcal D}_{H}({\mathfrak M})} } = (S^{{\mathcal D}_{N}({\mathfrak M})})_{{\mathfrak M}\cap S^{{\mathcal D}_{N}({\mathfrak M})} },$$ which shows 
(\ref{inclusion}). 

Let ${\mathfrak m}$ be any maximal ideal of $S^{N}$ and choose a maximal ideal $\tilde{\mathfrak M}$
of $S$ lying over ${\mathfrak m}$. Let $\sigma$ be an element of ${\mathcal D}_{H}({\mathfrak m})$.
As $\sigma(\tilde {\mathfrak M})$ is lying over ${\mathfrak m}$, we can choose
an element $\tau$ from $N$ in such a way that  $\tau (\tilde{\mathfrak M})
= \sigma(\tilde{\mathfrak M})$ (e.g.,  Chap. V, \cite{LR}).  Then by (\ref{inclusion}) we have $$\tau^{-1}\cdot \sigma
\in {\mathcal D}_{H}(\tilde{\mathfrak M}) \subseteq N,$$
which implies ${\mathcal D}_{H}({\mathfrak m})\subseteq N$.  

Next suppose that  $N$ contains
all  ${\mathcal I}_{H}({\mathfrak P})$  for  prime ideals  ${\mathfrak P}$
of $S$ of height $\leqq 2$ and $S^{H}$ is a locally complete intersection. 
We can similarly show (\ref{inclusion}) by Expose X of \cite{SGA2} 
 instead of Nagata's purity of branch loci.  The assertion follows similarly  from 
   the last paragraph as above.  \qed


\begin{theorem}\label{Chevalley-Serre} Suppose that $G^{0}$ is an algebraic torus. 
Let $(X, G)$ be an effective  regular action 
of   $G$ 
 on an affine Krull $K$-scheme
$X = {\rm Spec} (R)$ such that 
${\mathcal Q}(R^{G^{0}}) = {\mathcal Q}(R)^{G^{0}}$.   Suppose that
$R^{G}$ is regular. 
\begin{itemize}
\item[(i)]  If  $G = Z_{G}(G^{0})$ and there exists a maximal ideal ${\mathfrak M}_{0}$ of $R$
invariant under the action of $G$  such that 
$R/{\mathfrak M}_{0} \cong K$, 
then 
$G\vert_{R^{G^{0}}} = {\mathfrak R}(R, G)\vert_{R^{G^{0}}}.$
\item[(ii)] If $G = Z_{G}(G^{0})$ and $R$ is finitely generated over $K$ as a $K$-algebra
(i.e., $X$ is an affine normal variety), then
${\mathcal D}_{G}({\mathfrak m})\vert _{R^{N}}$
is trivial for any maximal ideal ${\mathfrak m}$ of $R^{N}$,
where $N=  G^{0}\cdot {\mathfrak R}(R, G)$. 
\item[(iii)] If $Z_{G}(G^{0}) \supseteq {\mathcal D}_{G}({\mathfrak m}_{0})$ 
for  a maximal ideal ${\mathfrak m}_{0}$ of $R^{G^{0}}$ such that $R^{G^{0}}/{\mathfrak m}_{0} 
\cong K$, then ${\mathcal D}_{G}({\mathfrak m}_{0}) \vert_{R^{G^{0}}} = 
{\mathfrak R}(R, {\mathcal D}_{G}({\mathfrak m}_{0}))\vert_{R^{G^{0}}}.$
\end{itemize}

\end{theorem}

\proof  We apply Proposition \ref{generalSerre}  to $S = R^{G^{0}}$,
 $H = G/{\rm Ker} (G \to {\rm Aut} (R^{G^{0}}))$ and  $N =
 {\rm Im} ({\mathfrak R} (R^{G^{0}}, G) \to H)$. By Theorem \ref{main} we must have
 $N = {\rm Im} ({\mathfrak R} (R, G) \to H)$. 
 The  action of 
 ${\mathcal D}_{H} ({\mathfrak M})$ 
is trivial on  $S/{\mathfrak M}$  for  a
maximal ideal ${\mathfrak M}$ of $S$ such that $S/{\mathfrak M}\cong K$.
 Thus the assertion in  {\it (ii)} follows from  Proposition \ref{generalSerre}. 
In the case where {\it (i)}    by Proposition \ref{generalSerre} again, 
we see  $H ={\mathcal D}_{H}({\mathfrak M}_{0}\cap S^{N})\subseteq N$,
which shows the assertion. 

{\it (iii):} As ${\mathcal D}_{G}({\mathfrak m}_{0})\vert_{R^{G^{0}}}
={\mathcal D}_{G/G^{0}}({\mathfrak m}_{0})\vert_{R^{G^{0}}}$,  by the proof of Proposition \ref{generalSerre}, we see that
$(R^{G^{0}})_{\mathfrak m_{0}\cap R^{{\mathcal D}_{G}({\mathfrak m}_{0})}}$
is a regular local ring. Applying Proposition \ref{generalSerre} to $S = (R^{G^{0}})_{{\mathfrak m}_{0}}$ 
acted by the finite group  ${\mathcal D}_{G}({\mathfrak m}_{0})\vert_{R^{G^{0}}}$, we must have
$${\mathcal D}_{G}({\mathfrak m}_{0})\vert_{R^{G^{0}}} ={\mathcal D}_{G}({\mathfrak m}_{0})\vert _{S}= {\mathfrak R}(S, {\mathcal D}_{G}({\mathfrak m}_{0}))\vert_{S} = 
{\mathfrak R}(R^{G^{0}}, {\mathcal D}_{G}({\mathfrak m}_{0}))\vert_{R^{G^{0}}}.$$
Then the assertion follows from Theorem \ref{main} for the action of ${\mathcal D}_{G}({\mathfrak m}_{0})$ on $R$, because $({\mathcal D}_{G}({\mathfrak m}_{0}))^{0} = G^{0}$
and $Z_{{\mathcal D}_{G}({\mathfrak m}_{0})}(({\mathcal D}_{G}({\mathfrak m}_{0}))^{0}) = {\mathcal D}_{G}({\mathfrak m}_{0})$. 
\qed

The assertion {\it (i)} or {\it (ii)} is  a generalization of the latter half
of the Chevalley-Serre Theorem and the assertion {\it (iii)} can be regarded as a
generalization of 
the Steinberg fixed point theorem (cf. \cite{Lehrer}).  
In fact Theorem \ref{Chevalley-Serre} is new,  even  if  
$R$ is a $K$-algebra of polynomial functions on $V$ of a finite dimensional
representation $G \to GL(V)$.


\medskip
       
\begin{thebibliography}{99}\small
       
           \bibitem{BourbakiC} N. Bourbaki,
             \newblock {Commutative Algebra : 
Chapters 1-7 
             (Elements of Mathematics)}
             \newblock Springer-Verlag, Berlin-Heidelberg-New York,  1989. 
             \bibitem{Bourbaki} N. Bourbaki,
             \newblock {Groupes et Alg\`ebres de Lie : Chapitres 4 \`a 6 
             (\'El\'ements de Math\'ematique)}
             \newblock Springer-Verlag, Berlin-Heidelberg-New York,  2006. 
               

         
            \bibitem{Kac} Victor. G. Kac; Kei.-Ichi. Watanabe,
            \newblock{Finite linear groups whose ring of invariants is a complete intersection,}
            \newblock{ Bull. Amer. Math. Soc. 6 (1982), 221-223.}
       \bibitem{Gr} Alexander Grothendieck; Mich\`ele, Raynaud, 
    \newblock {Rev\^etements \'Etales et Groupe Fondamental (SGA1), 
       Seminaire de geometrie algebrique du Bois-Marie 1960-61,}
      \newblock { Documents Math\'ematiques (Paris) 3, Soci\'et\'e Math\'ematique de France,  2003.}
   \bibitem{SGA2} Alexander Grothendieck; Mich\`ele, Raynaud, 
   \newblock{Cohomologie locale des faisccaux coh\'erents et
   Th\'eor/`emes de Lefschetz locaux et globaux (SGA2),}
   \newblock{North-Holland Publ. Company, Amsterdam 1968.}
   \bibitem{Lehrer} Gustav I. Lehrer, 
   \newblock{A new proof of Steinbergfs fixed-point theorem,}
   \newblock{Int. Math.  Res.  Notices  28 (2004), 1407-1411.}
        \bibitem{Magid} Andy Magid, 
             \newblock {Finite generation of class groups of rings
             of invariants,}
             \newblock Proc. of Amer. Math. Soc. 60 (1976), 45-48.
                       \bibitem{LR} Masayoshi Nagata,
             \newblock {Local Rings,}
             \newblock Robert E. Krieger Publ. Co., New York,  1975.
             
        \bibitem{Nak2} Haruhisa  Nakajima, 
             \newblock {Reduced ramification indices of quotient
             morphisms under torus actions,}
             \newblock J. Algebra 242  (2001), 536-549.
       \bibitem{Nak3} Haruhisa Nakajima,   \newblock {Divisorial free modules of relative invariants on Krull domains,} \newblock J. Algebra 292 (2005), 540-565.
             \bibitem{Nak5} Haruhisa Nakajima, 
       \newblock {Reductivities and finiteness of pseudo-reflections of algebraic groups
       and homogeneous fiber bundles,} \newblock {J. Pure and Appl. Algebra 217 (2013), 1548-1562.}
       \bibitem{Nak6} Haruhisa Nakajima, 
       {Valuative characterizations of central extensions of algebraic tori on Krull domains,  
       arXiv:1707.06008[math.GR],} to appear. 
        \bibitem{Serre} Jean-Pierre. Serre,
             \newblock{Groupes finis d'automorphismes d'anneaux locaux r\'eguliers,}
             \newblock{Colloque d'Alg\`ebre (Paris, 1967), Exp. 8,  11 pp,}
             \newblock{ Secr\'etariat math\'ematique, Paris,  1968.} 
             \bibitem{Stanley} Richard P. Stanley,
             \newblock{Relative invariants of finite groups generated by pseudoreflections,}
             \newblock{J. Algebra 49 (1977), 134-148.}
        \bibitem{Po} Vladimir L. Popov; Tonny A. Springer,
        \newblock{Algebraic Geometry IV,}
         \newblock{Encyclopaedia of Mathematical Sciences 55,}
        \newblock{Springer-Verlag, Berlin-Heidelberg-New York, 1994.}
     \end{thebibliography}
\end{document}